\documentclass{tufte-handout}

\title{A Very Elementary Introduction to Sheaves}
\author{Mark Agrios}
\date{\today}

\usepackage{amsmath}
\usepackage{embedall}
\usepackage{tikz-cd}
\usepackage{booktabs}
\usepackage{units}
\usepackage{fancyvrb}
\fvset{fontsize=\normalsize}
\usepackage{multicol}
\usepackage{lipsum}
\usepackage{graphicx}
\graphicspath{{img/}}
\usepackage{placeins}
\usepackage{hyperref}
\usepackage{cleveref}
\usepackage{titlesec}

\setkeys{Gin}{width=\linewidth,totalheight=\textheight,keepaspectratio}

\newcommand{\F}{\mathcal{F}}
\newcommand{\R}{\mathbb{R}}

{\vspace{1em}\begin{center}\Large\itshape}{\end{center}\vspace{1em}}

\hypersetup{
    colorlinks=true, 
    linktocpage=true,     
    linkcolor=Maroon,  
    urlcolor=Maroon,
}
\begin{document}
\maketitle

\begin{abstract}
\noindent Here, we give a simple and approachable explanation to the mathematical objects called sheaves. This paper presents tangible and concrete examples so that readers will be able to further explore these concepts on their own in more abstraction and application.
\end{abstract}


\marginnote[0pt]{There are things called \emph{pre}-sheaves which are a less specific object used to define a true sheaf. It just isn't helpful to bring up that distinction in an introductory paper. So sometimes when the word "sheaf" is used it might actually mean presheaf.} 

\subsection{About This Paper}
This paper is a very \emph{non}-rigorous, loose, and extremely basic introduction to sheaves. This is meant to be a a guide to gaining intuition about sheaves, what they look like, and how they work, so that after reading this paper, someone can jump into the extremely abstract definitions and examples seen in textbooks with at least some idea of what is going on. Most of this material is inspired and built from the work of Dr. Michael Robinson\cite[-40pt]{robinsonTSP}\cite[-18pt]{robinsonPOSET}\cite{robinsonSheaf}, and that of Dr. Robert Ghrist and Dr. Jakob Hansen\cite{hansengentle}\cite{hansen2019learning}, as well as Dr. Justin Curry's PhD thesis%
\cite{Curry},
who are some of the only applied sheaf theorists out there and they do an amazing job of explaining sheaves in a concrete way through their research. The rest of this paper is populated by mathematical definitions found in textbooks that I have stretched from two lines into multiple pages, as well as some analogies for thinking of sheaves I have thought of myself. This paper only assumes knowledge of basic linear algebra, basic group theory, and the very fundamentals of topology. If there is anything in the setup that you do not understand it is probably a quick Wikipedia search away. I hope this paper provides insight, intuition, and helpful examples of why sheaves are such powerful tools in both math and science.

\section{A Quick Overview}
A sheaf can be thought of as a way to "enhance" some kind of mathematical object that you are given. Think about it like the mathematical object is a plot of land and a sheaf is like a garden on top of it. On one square of dirt, you could organize a garden in many different ways, and once you have built your garden then the act of harvesting it becomes its own endeavor that is dependent on how you organized it%
\sidenote{These two parts --- laying out a garden and harvesting from the garden --- are important notions to keep in mind. Sheaves have two parts, the framework you build and what you can do with within this framework.}%
. The main goal from this paper is to drive home the point that sheaves take a single, rigid, mathematical object and build on top of it a framework that is both flexible and powerful but still remains faithful to the underlying structure of the object in question, just like how you can do a lot more with a garden than a plot of dirt, but the plot of dirt is essential to make a garden. 

To build a sheaf on a mathematical object, you need to understand what pieces make up that object and how those pieces are arranged. Graphs are made of edges that connected nodes, topological spaces are made up of open sets. A common idea when studying sheaves is that of \emph{direction}. When we say that a mathematical object has a certain "structure" that the sheaf must respect, we usually mean the direction of how the pieces fit inside one another. Edges on a graph are often thought of as the intersection of the nodes it connects, so there is a partial order on a graph where $e_{ij} < v_i,v_j$ for an edge $e_{ij}$ connecting nodes $v_i$ and $v_j$. In a topological space, there is a natural partial order of the open sets. For an open set $U$, and an open set $V$ contained in $U$, the order is $V \subset U$. It is this sense of direction where you will often hear people talk about how a sheaf allows us to attach both \emph{local} and \emph{global} \emph{data}%
\sidenote{Our mathematical object is the plot of dirt, our sheaf is the garden we build on top of it. Think of "data" as the vegetables grown in this garden. Maybe think of  "local data" as an apple and "global data" being all the apples from one tree. This will make more sense in a bit, don't worry about it too much.}
to features of an object. 

\begin{figure}[h!]
    \centering
    \includegraphics{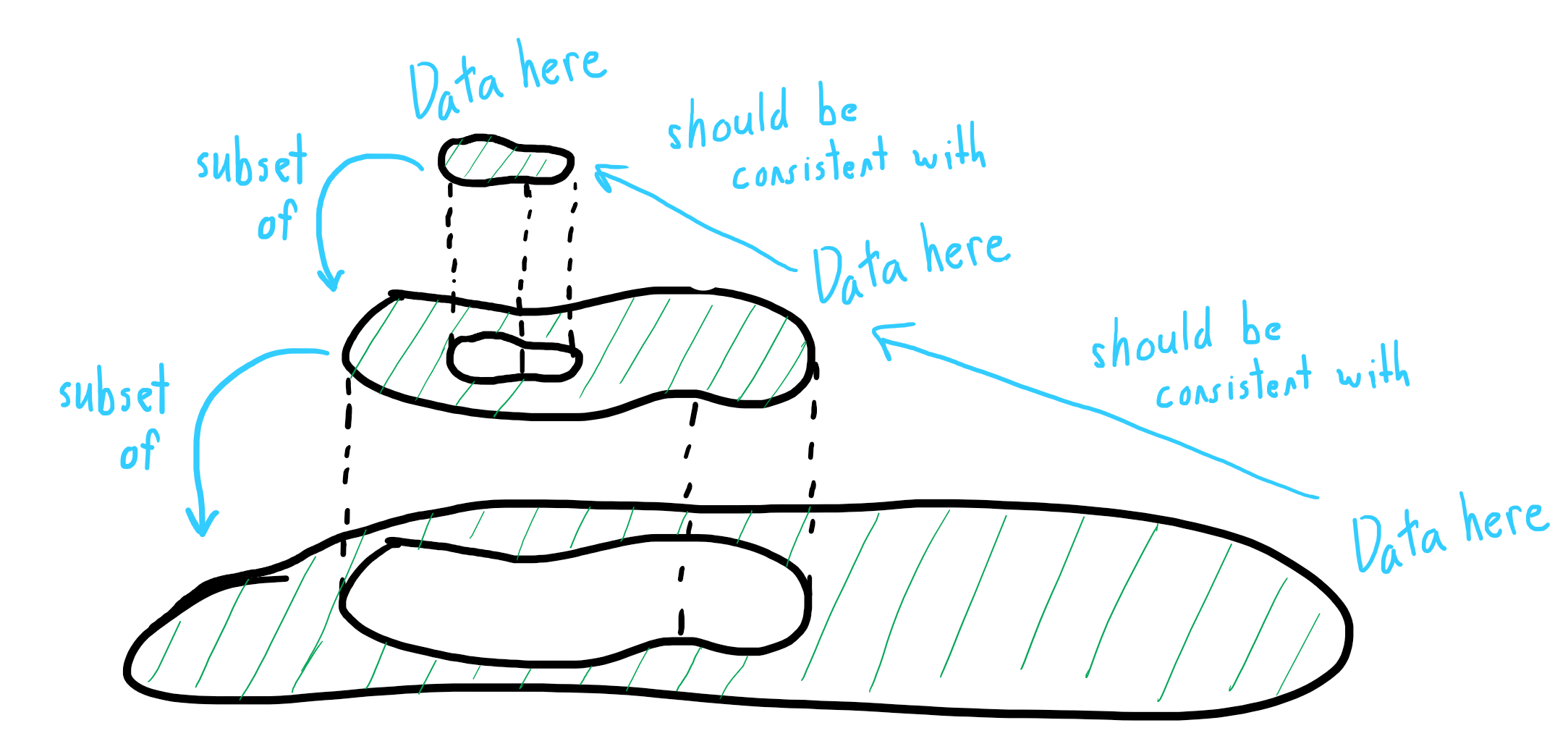}
    \caption{The mathematical object on which our sheaf is built has a structure that is described by how pieces fit together locally to globally. A sheaf organizes data attached to these pieces in a way that creates consistency from global to local.}
    \label{fig:refine}
\end{figure}

If you held this paper up side by side to another more mathematical or abstract introduction to sheaves, this might be the part where it gets confusing. Sheaves take this order (local $\rightarrow$ global) and reverse it. Data attached to a larger part of the object must be consistent when restricted to a smaller part of the of the object (global $\rightarrow$ local). Think about it like reading a novel. If there is a sub-plot that spans one chapter then the actions and dialogue of the characters at a certain paragraph in the chapter must make sense in the context of that chapter. Similarly for a sentence within that paragraph. This notion is illustrated in \cref{fig:refine} and is often referred to as \emph{refinement} or \emph{increasing resolution}. 

In this paper we will first take a look at sheaves on graphs, because they are finite and discrete and then later build intuition about how to build sheaves on topological spaces.

\section{Now an Example}
Here, we are going to introduce an example of a sheaf before the strict definition of a sheaf. This is done because the mathematical definition of a sheaf is quite cumbersome and confusing, and it helps to have a picture in mind when reading the definition. So let's start with the picture.

\begin{marginfigure}%
  \includegraphics[width=\linewidth]{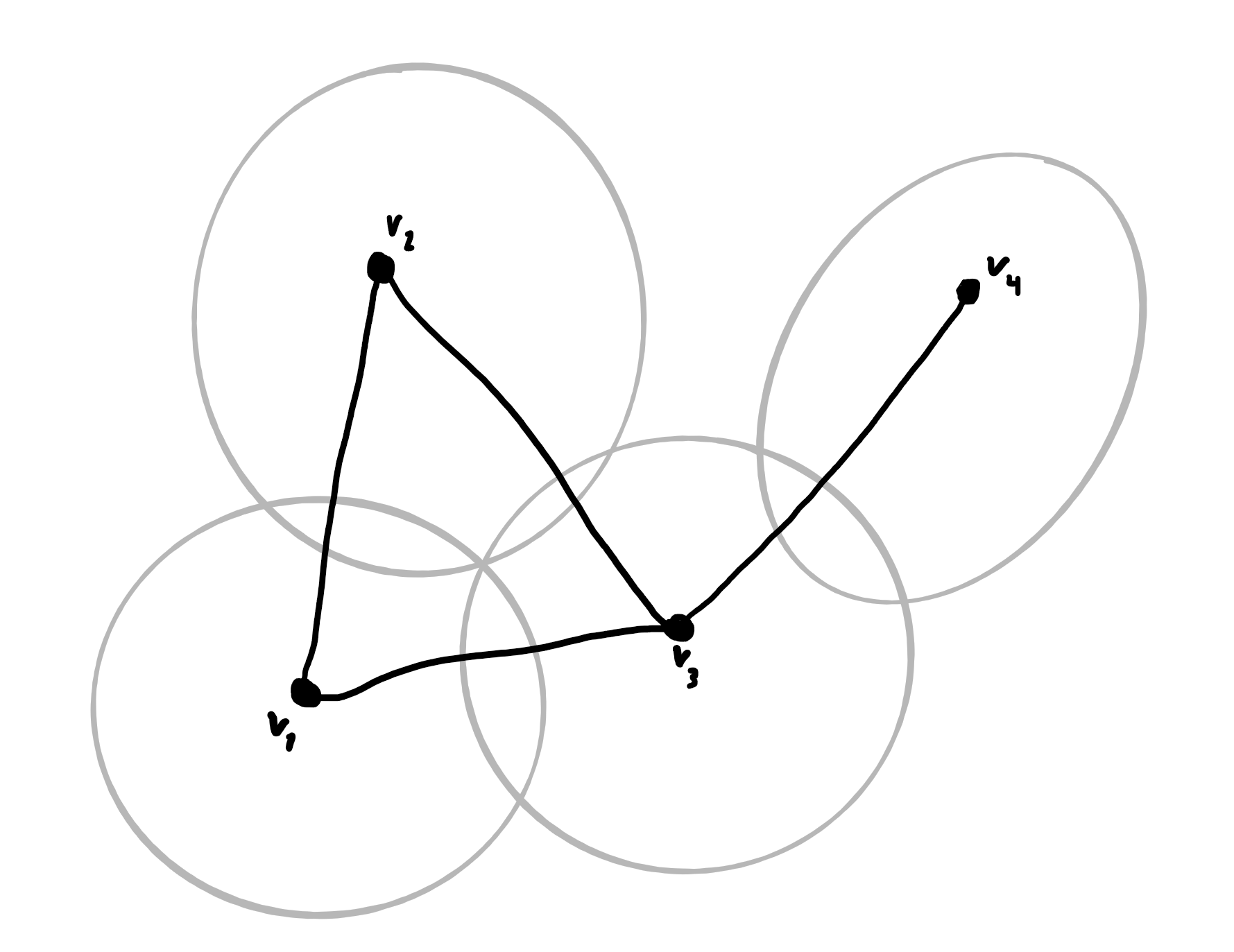}
  \caption{We can think of an edge as the intersection of two nodes. This lets us describe our graph in the language of set theory which our sheaf will appreciate greatly.}
  \label{fig:setedge}
\end{marginfigure}

\begin{figure}[h!]
    \centering
    \includegraphics{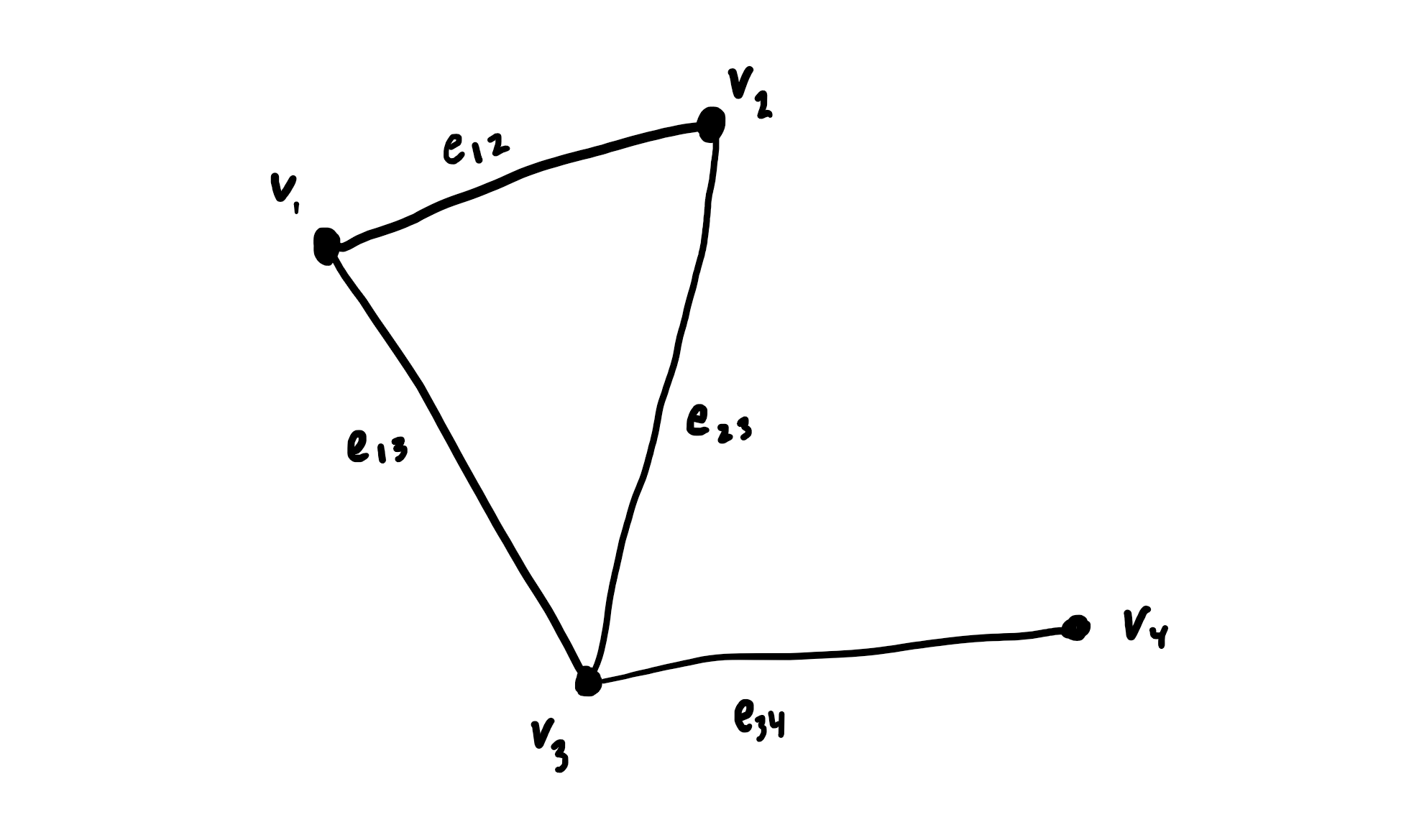}
    \caption{A graph with nodes $v_1,...v_4$ and edges $e_{12},...,e_{34}$, and so the \emph{components} of this graph are $\{v_1,...v_4,e_{12},...e_{34}\}$, where each $e_{ij} \subseteq v_i,v_j$.}
    \label{fig:graph}
\end{figure}

Suppose we have a graph. The \emph{components}, as we will call them, of this graph are the collection of nodes and edges we can see labeled in \cref{fig:graph}. These are the building blocks of a graph and will be foundation on which our sheaf lives. Let's build a sheaf, call it $\F$, on this graph.

\subsection{The First Mechanism}

A sheaf has two mechanisms. The first mechanism takes each component and assigns it a new space called the \emph{stalk}. So for a node $v_i$, $\F$ assigns it to, say, the vector space of real numbers $\R$. In words, we say "the stalk over $v_i$ is $\R$" and we write this as 
\[
    \F(v_i) = \R
\]
This association of a component of the graph to new space has a lot of freedom. A component can be assigned to just about anything as long as all the components play nicely with each other (as we will see later). So we could have something totally random like
\begin{align*}
    \F(v_1) &= \R^2 \\
    \F(v_2) &= \R^{100} \\
    \F(e_{12}) &= \R^{34} \\
    \vdots
\end{align*}  
Or, $\F$ could assign every node and edge to $\R$, or assign each $v_i$ to $\R^{deg(v_i)}$ and all the edges to $\R$ where $deg(v_i)$ is the number of edges that come out of node $v_i$. Really all that matters is every component needs to be mapped to some kind of space.

\subsection{The Second Mechanism}

Let's now address this notion of requiring the stalks to "play nicely together". This simply means that all the stalks need to be similar enough so that we can define functions that go between them. In the example above, all of our stalks are a vector space $\R^n$ for some integer $n$. This leads us to the second mechanism of a sheaf. 

For a sheaf $\F$ of a graph, and an edge $e_{ij}$ connecting nodes $v_i$ and $v_j$, $\F$ has two functions that go from the stalk over each node to the stalk over the edge, namely
\begin{align*}
    \F_{v_i:e_{ij}}: \F(v_i) \rightarrow \F(e_{ij}) \\
    \F_{v_j:e_{ij}}: \F(v_j) \rightarrow \F(e_{ij})
\end{align*}

\begin{figure}[h!]
    \centering
    \includegraphics{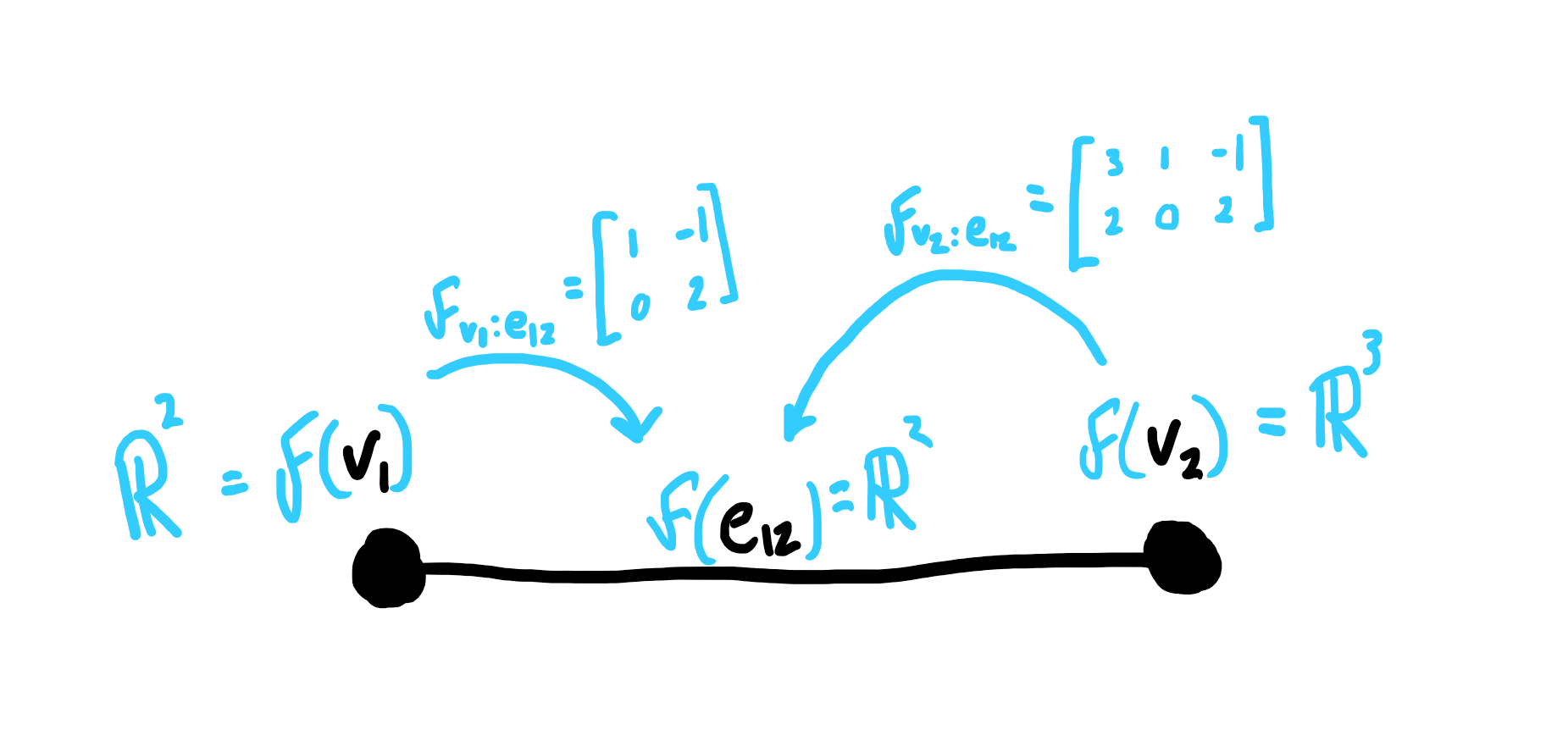}
    \caption{Part of a graph (black) with the corresponding parts of the sheaf (blue). Here we can see that each of the node stalks have restrictions maps into the same edge stalk. This is because the edge $e_{12}$ is thought of as the intersection of nodes $v_1$ and $v_2$.}
    \label{fig:sheafassignspace}
\end{figure}

And the sheaf has these maps for every edge and incident nodes. These maps are called the \emph{restriction maps}%
\sidenote{They are called restriction maps because they restrict some piece to a smaller part of that piece. Remember the global $\rightarrow$ local direction mentioned earlier; in the context of graphs, an edge is a subset of both the nodes that it connects.}
of the sheaf $\F$. Let's refer back to our graph and just look at one edge $e_{12}$ that connects $v_1$ to $v_2$. Suppose $\F$ has the following stalks over these components
\begin{align*}
    \F(v_1) &= \R^2 \\
    \F(v_2) &= \R^3 \\
    \F(e_{12}) &= \R^2
\end{align*}
Since each stalk is a vector space, it is natural to have our restriction maps be linear transformations described by matrices. It would be a little overwhelming to draw every stalk and restriction map over our graph, so \cref{fig:sheafassignspace} illustrates just one a piece of our graph with the relevant parts of $\F$, but there are stalks and restriction maps for all components.

And we have complete described our sheaf! It is a map from components of a mathematical structure to separate spaces along with functions that let you move between them. In words, we would say "$\F$ is a sheaf of \emph{vector spaces} on a \emph{graph}".


\marginnote[-4em]{In the literature, you will often come across sheaves of abelian groups on a topological space, or sheaves of rings on a topological space. This is because sheaves were originally invented within the context of algebraic topology and then later found to be very helpful in algebraic geometry. We will handle these types of sheaves in the next section.}%

\subsection{Using Our Sheaf}

So we have vector spaces, called stalks, attached to each component of our graph. In mathematics, we work in certain spaces, and manipulate elements of these to find interesting results. So let's just pick out some vectors from each stalk and see what happens. Again, we will just illustrate a piece of our graph but everything we say here applies to our entire graph and sheaf. 

\begin{figure}[h!]
    \centering
    \includegraphics{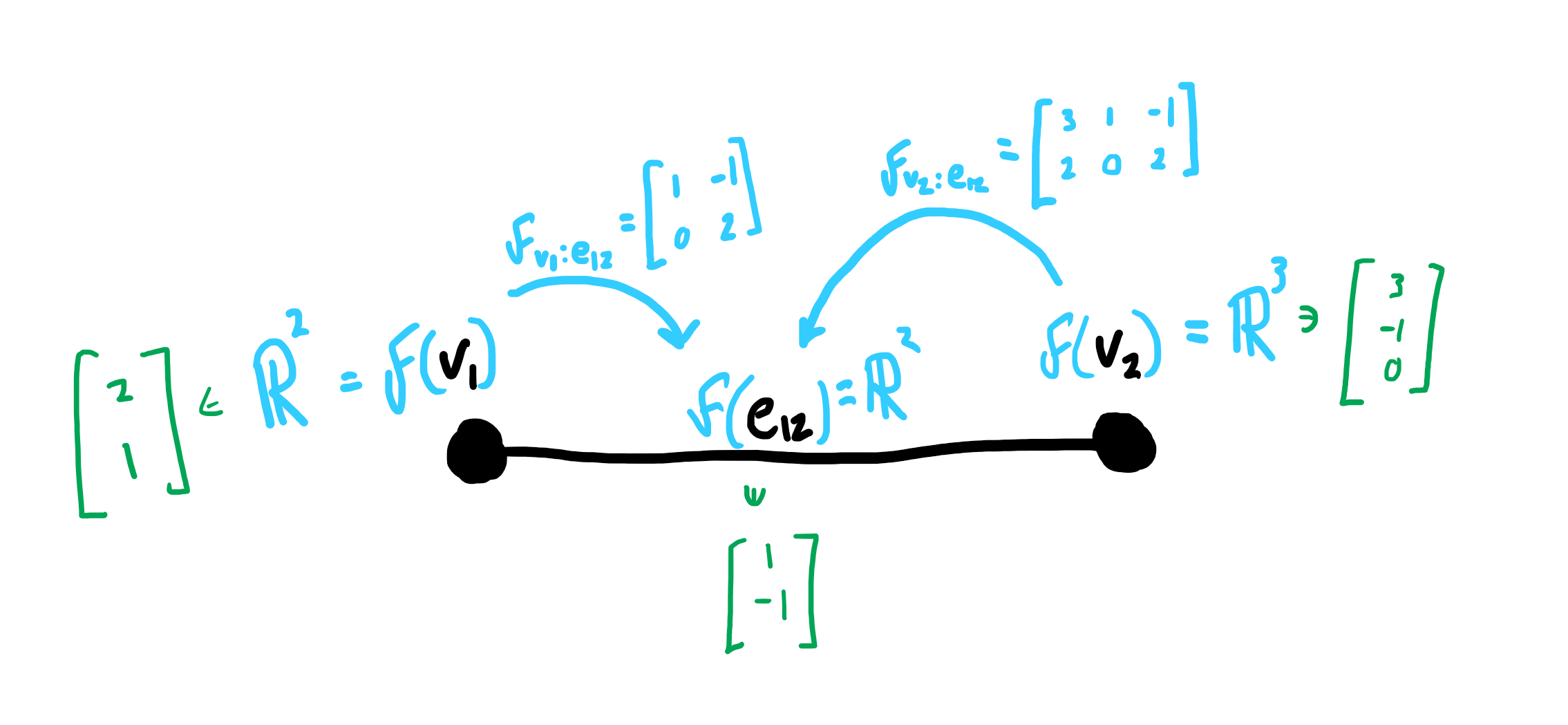}
    \caption{A section (green) of our sheaf $\F$ (blue) on our graph. Note that this section is totally random, we just arbitrarily chose vectors from each stalk.}
    \label{fig:section}
\end{figure}

A choice of element from each stalk of a sheaf (\cref{fig:section}) is called a \emph{section}. Zooming out a level of abstraction, if we think of a sheaf as one coherent space, then a section is one element of this space, just like vectors are elements of vector spaces and points are elements of topological spaces. Let's look at this section we just picked out.

Using the tools of our sheaf, it is natural to think about what happens when we move these vectors through their restriction maps. Consider the first one $\large[ \begin{smallmatrix} 2\\1 \end{smallmatrix}\large]$. 
\[
    \F_{v_1:e_{12}}\Bigg(\begin{bmatrix} 2\\1 \end{bmatrix}\Bigg) = %
    \begin{bmatrix}
        1 & -1 \\
        0 & -2
    \end{bmatrix}%
    \begin{bmatrix}
        2 \\
        1
    \end{bmatrix}=%
    \begin{bmatrix}
        1 \\
        -2
    \end{bmatrix} \neq%
    \begin{bmatrix}
        1 \\
        -1
    \end{bmatrix}    
\]
We can see that the choice of $\large[ \begin{smallmatrix} 2\\1 \end{smallmatrix}\large]$ in the stalk $\F(v_1)$ is not consistent with the choice of $\large[ \begin{smallmatrix} 1\\-1 \end{smallmatrix}\large]$ in the stalk $\F(e_{12})$ because  $\F_{v_1:e_{12}}(\large[ \begin{smallmatrix} 2\\1 \end{smallmatrix}\large]) \neq \large[ \begin{smallmatrix} 1\\-1 \end{smallmatrix}\large]$. Similarly for the choice of vector in $\F(v_2)$ as one can check.

\begin{figure}[h!]
    \centering
    \includegraphics{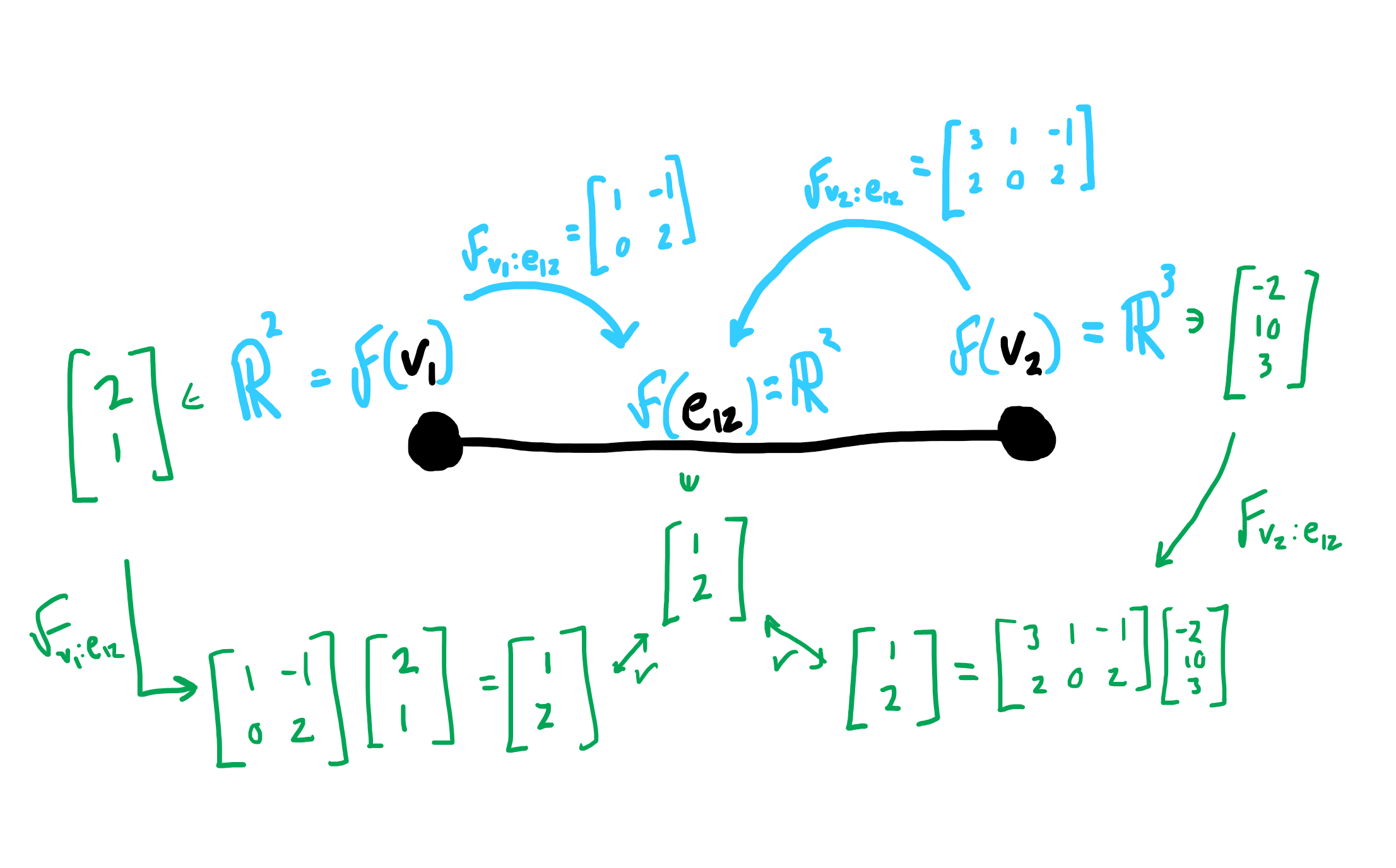}
    \caption{A piece of the graph (black) with the sheaf (blue) and the global section (green) illustrating how the vectors at each node align with the vector at the edge through the restriction maps.}
    \label{fig:globalsection}
\end{figure}

Now, wouldn't it be interesting if we were to choose a vector at each stalk of $\F$ such that for every edge $e_{ij}$, 
\[
\F_{v_i:e_{ij}}(x_i) = y_{ij} = \F_{v_j:e_{ij}}(x_j)
\]
For every $x_i \in \F(v_i), x_j \in \F(v_j) $ and $y_{ij} \in \F(e_{ij})$. That is, the choice of vector at every stalk is consistent with its image under its restriction maps. That would be called a \emph{global section}. The term global is used because we look at the data associated with the entire space and see how it remains consistent when restricting it to a more local piece. Again, not to draw the entire sheaf over our graph, a piece of a global section is illustrated in \cref{fig:globalsection}. Global sections are the important and interesting types of data that we are tasked to find within our sheaf. A sheaf can have infinitely many global sections or it can have none at all, or it could have eleven or just one. It all depends on how the sheaf is setup. In practice, many questions are answered by finding out what the set of all global sections of a sheaf, usually denoted $\Gamma(\F)$, looks like. Just like how for a group homomorphism, the kernel is where a lot of the interesting stuff is.

So this is an example of a sheaf of vector spaces on a graph. Most likely, you will not see sheaves explained this way, but having an example that is discrete and finite makes it a lot easier to grasp the basic ideas of what a sheaf is supposed to do. Most textbooks will introduce sheaves as sheaves of abelian groups on a continuous topological space, as you can see this is a continuous example with the underlying partial order (the open sets) being infinite%
\sidenote[][-30pt]{It wasn't explicitly mentioned earlier but in fact the \emph{only} structure you need to construct a sheaf is a partial order. Graphs and topological spaces have very obvious partial orders, but it is possible to make a sheaf on just a partial order.}. 

\section{An Example You'll See in a Textbook that We are Extending to Multiple Pages}

For a topological space $X$, a sheaf of abelian groups on $X$ assigns each open set of $X$ to an abelian group. Which abelian group you ask? Here, it will be the group of continuous functions from $U$ to $\R$. That is, for an open set $U \subseteq X$, the stalk of $U$ is the group
\[
    \F(U) := \{f: U \rightarrow R: f \text{ is continuous}\}
\]
along with point-wise addition as the group operation. And for an open set $V \subseteq U$ the restriction map $\F_{U:V}: \F(U) \rightarrow \F(V)$ is actually weirdly simple. For an element $f$ of $\F(U)$, which is a continuous function, $\F_{U:V}(f)$ is the restriction of $f$ from $U$ to $V$, often denoted $f|_V$. So $\F_{U:V}$ is a group homomorphism between $\F(U)$, the elements of which are continuous functions from $U$ to $\R$, and $\F(V)$, the elements of which are continuous functions from $V$ to $\R$.
\[
    \F_{U:V}: \F(U) \rightarrow \F(V)
\]
where $f \mapsto f|_V$. This seems almost unhelpful to have a function just be mapped to itself over a smaller domain, but we'll see why we construct a sheaf this way. 

First, let's say the topological space we want to build a sheaf on is $\R$. So each stalk on an open set is the group of continuous functions from a subset of $\R$ to $\R$. This makes it easy to draw. 

\begin{figure}[h!]
    \centering
    \includegraphics{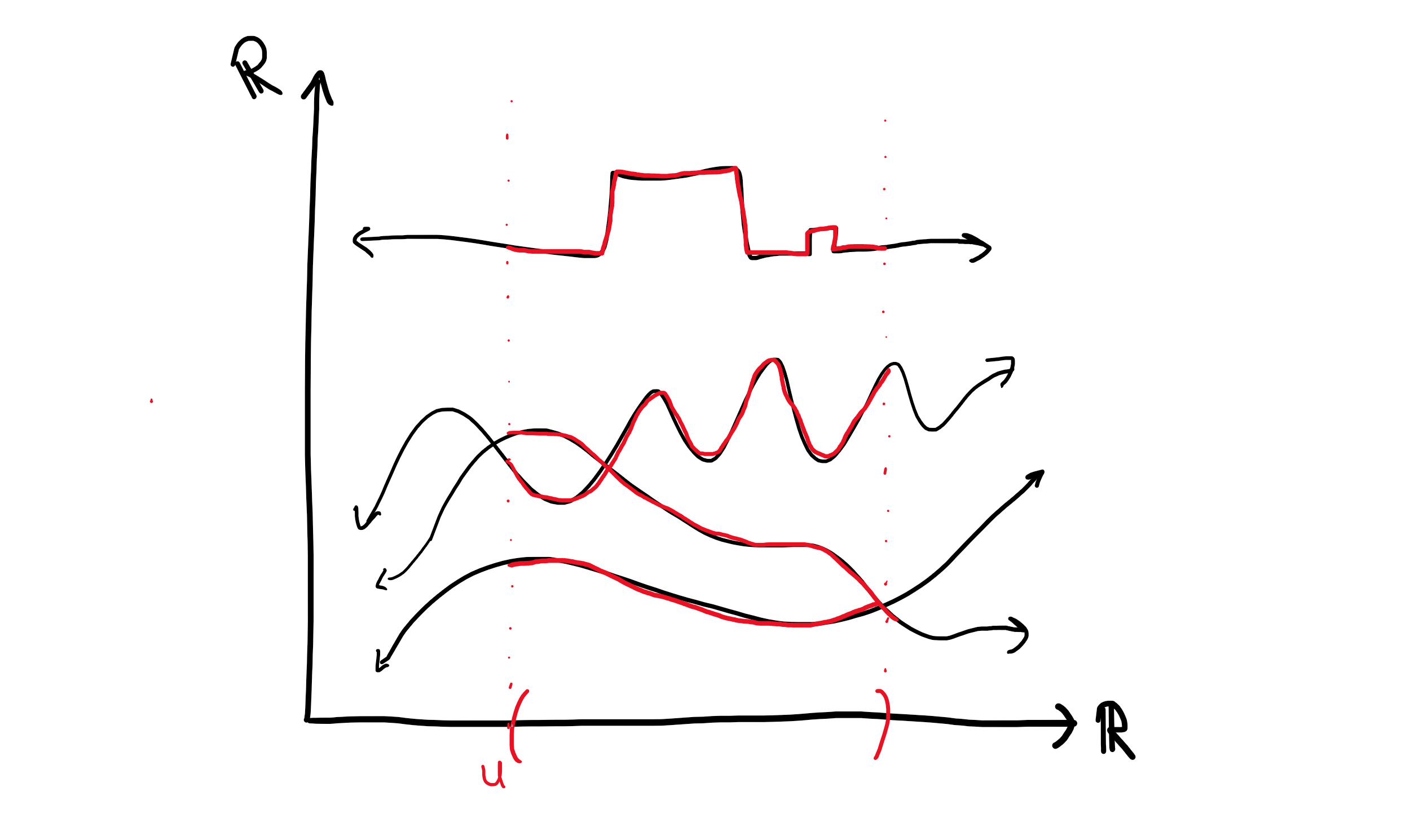}
    \caption{For an open set $U$ of $\R$, $\F(U)$ is only concerned with continuous functions from $U$ to $\R$, which means only the red part of these functions are elements of $\F(U)$. In fact, they don't even have to be defined outside of $U$.}
    \label{fig:contonR}
\end{figure}

\begin{figure}[h!]
    \centering
    \includegraphics{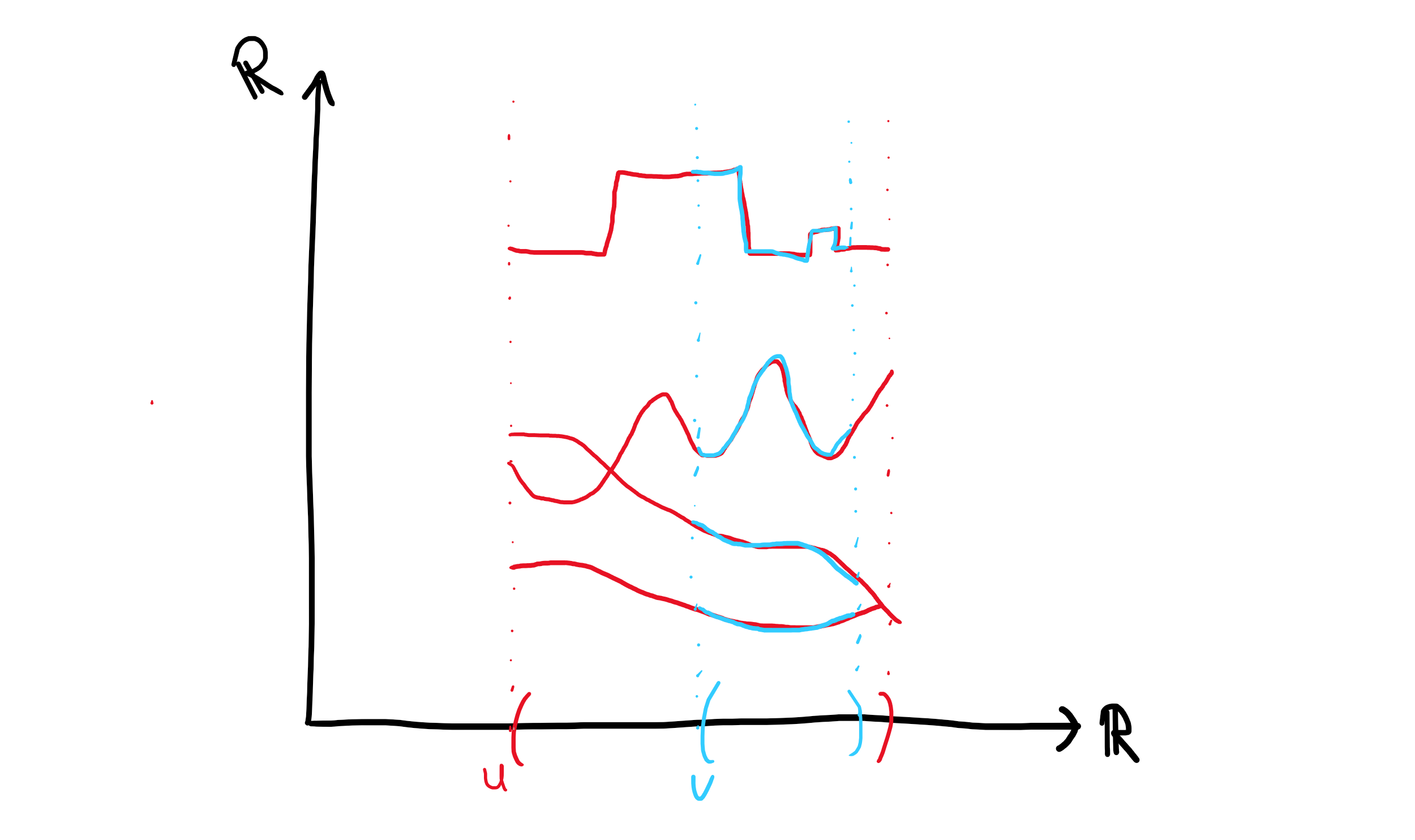}
    \caption{Each function in red is an element of $\F(U)$ while that function's image under $\F_{U:V}$ is shown in blue and is an element of $\F(V)$.}
    \label{fig:contonR_res}
\end{figure}

The red lines in \cref{fig:contonR} are some of the elements of $\F(U)$. For an open set $V \subseteq U$, the images of each $f \in \F(U)$ under $\F_{U:V}$ are shown in \cref{fig:contonR_res}.

Remember, we are working algebraically, so there is no analytic manipulation of the functions here, so the points inside $U$ but outside $V$ aren't mapped to the boundary of $V$ or anything like that. Each function is its own element in a group and it is being mapped to another element in a different group, it just happens to be that the element and its image under this group homomorphism look very similar to us. 

\marginnote[0pt]{Maurice Auslander said \emph{'Sheaf theory is the subject in which you do topology horizontally and algebra vertically.'} }

\begin{figure}[h!]
    \centering
    \includegraphics{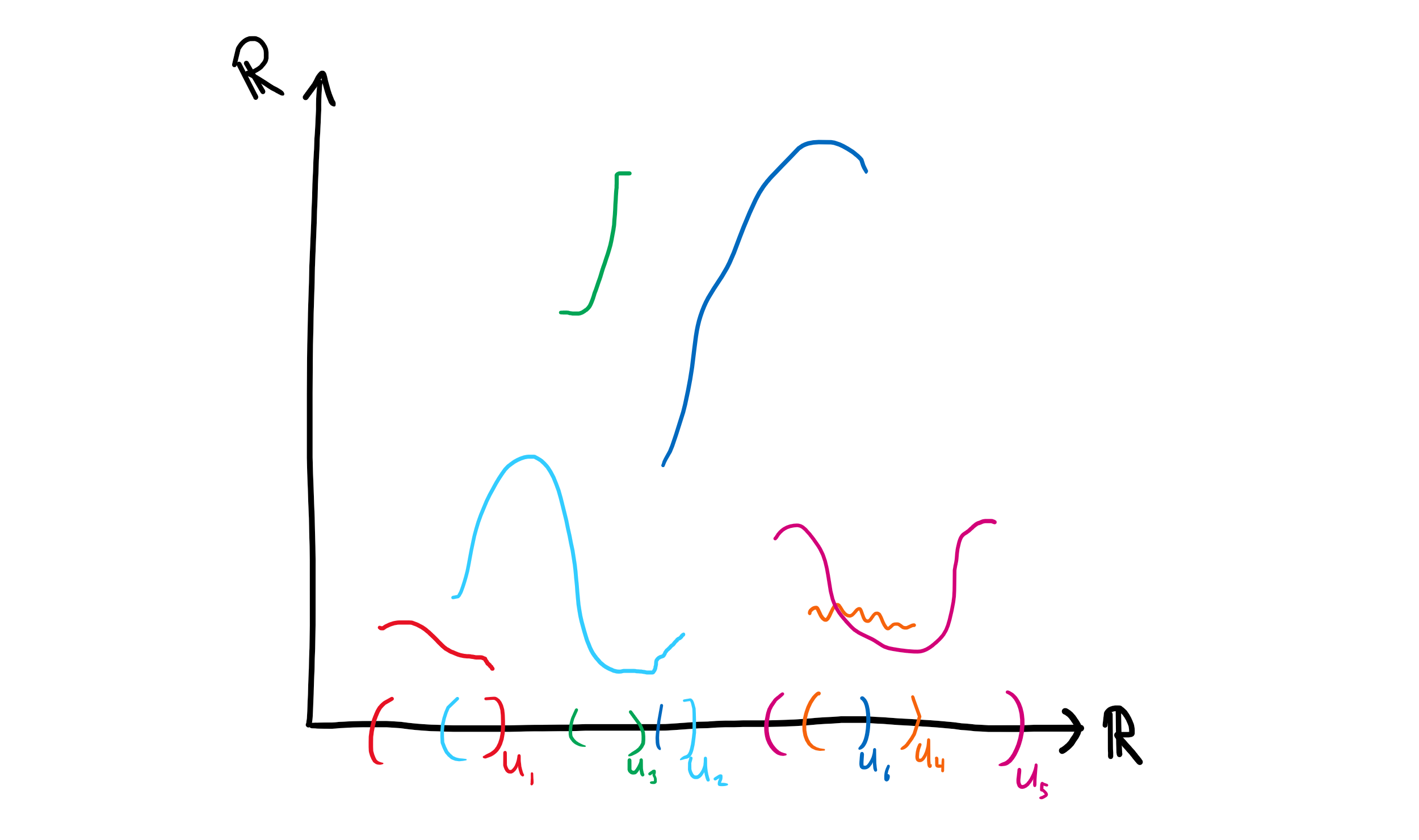}
    \caption{Here, we illustrate some pieces of a section of our sheaf, since we can't draw all uncountably many of them.}
    \label{fig:sectiononR}
\end{figure}

So now, what would a section look like? As we know, a section is a choice of element of each stalk. For our topological space, there are are infinitely many open sets (uncountably many, in fact) and hence infinitely many stalks. So a section is a choice of one function from $\F(U)$ for every open set $U$ of $\R$. That could look something like \cref{fig:sectiononR}.

\FloatBarrier

But, just like in our first example, this is not particularly helpful. And as we can see for the green line, call it $f_3$, and the light blue line, call it $f_2$, $\F_{U_2:U_3}(f_2) \neq f_3$, which is one of the many inconsistencies in this section. So again, the natural next step is to ask what a global section would look like. It would have to agree on every piece of the topological space. More symbolically, it must fulfill these criteria.
\begin{enumerate}
    \item For open sets $U \subseteq V \subseteq W$, $\F_{V:U} \circ \F_{W:V} = \F_{W:U}$. This means that data defined at some larger piece must be consistent with every refinement of that piece (see \cref{fig:refine}). Remember the chapter $\rightarrow$ paragraph $\rightarrow$ sentence analogy.
    
    \begin{marginfigure}%
      \includegraphics[width=\linewidth]{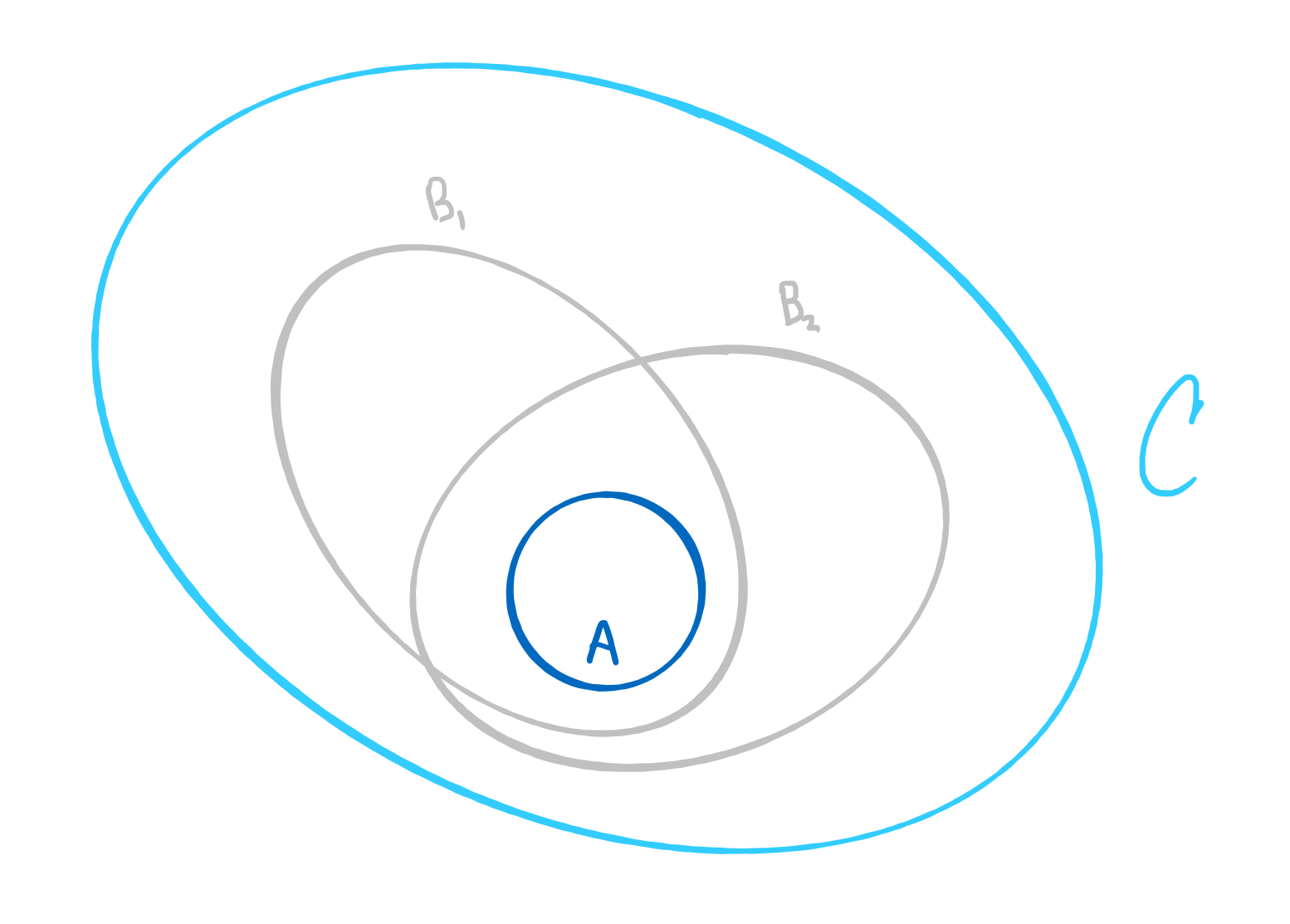}
      \caption{An illustration of two different refinements, $A \subseteq B_1 \subseteq C$ and $A \subseteq B_2 \subseteq C$.}
      \label{fig:tworefine}
    \end{marginfigure}
    
    \item For two orders of open sets, say $A \subseteq B_1 \subseteq C$ and $A \subseteq B_2 \subseteq C$, the composition of the restriction maps must be equal. That is, $\F_{B_1:A} \circ \F_{C:B_1} = \F_{B_2:A} \circ \F_{C:B_2}$. Meaning that if you want to refine to a certain open set, it does not matter how you get there.
    
    \item For open sets $U_1,U_2$ such that $U_1 \cap U_2 := U_{12} \neq \emptyset$, and for all elements $f_1 \in \F(U_1), f_2 \in \F(U_2)$, we have $\F_{U_1:U_{12}}(f_1) = \F_{U_2:U_{12}}(f_2)$ for all pairs of open sets with nonempty intersection. As we demonstrated in our graph example.

\end{enumerate}

\begin{figure}[h!]
    \centering
    \includegraphics{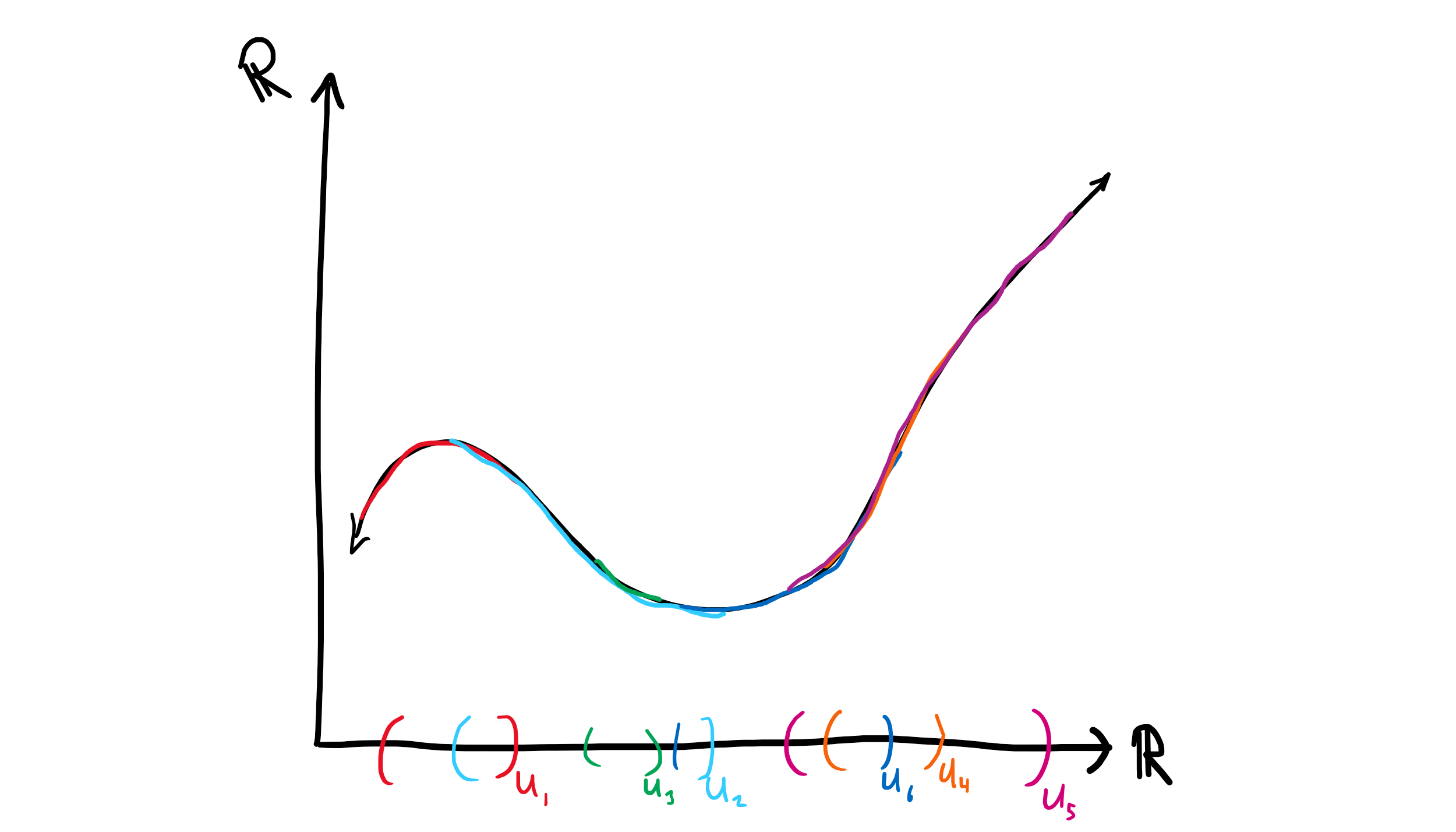}
    \caption{One possible global section of our sheaf. Notice how for every $f \in \F(U)$ and $g \in \F(V)$ where $V \subseteq U$ the consistency of $\F_{U:V}(f) = g$ holds.}
    \label{fig:globalsectiononR}
\end{figure}

Which might look something like \cref{fig:globalsectiononR}

\FloatBarrier

This is amazing! All we did was stipulate what must happen locally and from those rules we have created something with global structure: a continuous function from the entirety of $\R$ to $\R$! Remember how we said earlier that the space of global sections is usually what we are most interested in? With this sheaf of abelian groups we have built on $\R$, we have determined that the space of global sections of $\F$, remember we call this $\Gamma(\F)$, is all the continuous functions from $\R$ to $\R$. This could have been done with $k$-differentiable functions, smooth functions, constant functions, or a million other things, and we didn't have to build our sheaf on $\R$, we just did that so it was easier to draw, but we could build a sheaf on any topological space, or even just a partial order. There is an incredible amount of flexibility of what you can build a sheaf on and what that sheaf can look like. If you go on into more algebraic geometry, algebraic topology, differential geometry, or just math in general, you'll find sheaves make appearances in many different places with totally different applications. 

\section{Further Reading}

\begin{marginfigure}%
  \includegraphics[width=\linewidth]{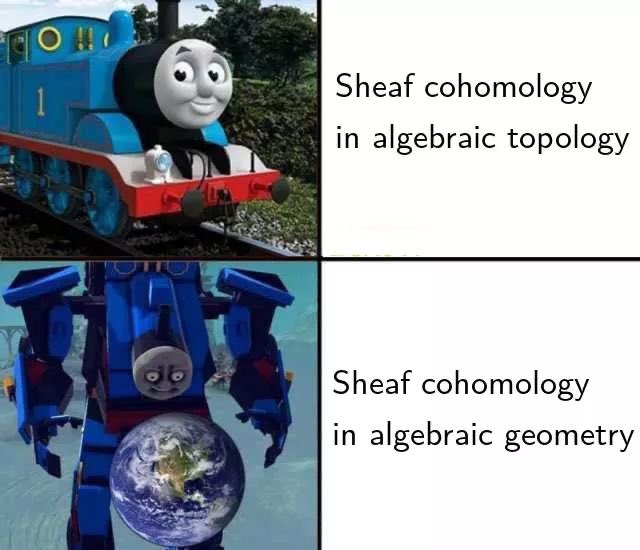}
  \caption{A visual warning to start your sheaf journey within algebraic topology before algebraic geometry. \emph{Source:} Reddit.}
  \label{fig:meme}
\end{marginfigure}

I encourage you to read the Wikipedia article%
\sidenote[][10pt]{\url{https://en.wikipedia.org/wiki/Sheaf_(mathematics)}}
on sheaves and see how far you can get. If you want to dive deeper into how sheaves are used outside of pure mathematics, I highly recommend looking more into the work done by Dr. Michael Robinson, Dr. Robert Ghrist, and Dr. Jakob Hansen. And even if you want to explore sheaves for purely mathematical reasons, check out Dr. Justin Curry's thesis, \emph{Sheaves, cosheaves, and applications}, although, I would still recommend the applied sheaf literature because they will present sheaves with thorough explanations and concrete examples that you can abstract later down the road. 

\bibliography{sources}
\bibliographystyle{plainnat}

\end{document}